# Impulsive fishery resource transporting strategies based on an open-ended stochastic growth model having a latent variable

Running title: Impulsive fishery resource transport

**Authors**

Hidekazu Yoshioka[1, 2, *], Tomomi Tanaka[2], Futoshi Aranishi[1, 2], Motoh Tsujimura[3], Yumi Yoshioka[1]

[1]Graduate School of Natural Science and Technology, Shimane University, 1060 Nishikawatsu, Matsue, 690-8504, Japan

[2]Fisheries Ecosystem Project Center, Shimane University, 1060 Nishikawatsu, Matsue, 690-8504, Japan

[3]Graduate School of Commerce, Doshisha University, Karasuma-Higashi-iru, Imadegawa-dori, Kamigyo-ku, Kyoto, 602-8580, Japan

*Corresponding author: yoshih@life.shimane-u.ac.jp

**Acknowledgements**

JSPS Research Grant 19H03073, Kurita Water and Environment Foundation Grant No. 19B018 and No. 20K004, Grant for Environmental Research Projects from the Sumitomo Foundation No. 203160, and grants from MLIT Japan for ecological survey of a life history of the landlocked *Plecoglossus altivelis altivelis* and management of seaweed in Lake Shinji support this research. The authors thank Dr. Masahiro Horinouchi in Shimane University for his valuable comments and suggestions on the manuscript.

**Abstract**

In inland fisheries, transporting fishery resource individuals from a habitat to spatially apart habitat(s) has recently been considered for fisheries stock management in the natural environment. However, its mathematical optimization, especially finding when and how much of the population should be transported, is still a fundamental unresolved issue. We propose a new impulse control framework to tackle this issue based on a simple but new stochastic growth model of individual fishes. The novel growth model governing individuals' body weights uses a Wright-Fisher model as a latent driver to reproduce plausible growth dynamics. The optimization problem is formulated as an impulse control problem of a cost-benefit functional constrained by a degenerate parabolic Fokker-Planck equation of the stochastic growth dynamics. Because the growth dynamics have an observable variable and an unobservable variable (a variable difficult or impossible to observe), we consider both full-information and partial-information cases. The latter is more involved but more realistic because of not explicitly using the unobservable variable in designing the controls. In both cases, resolving an optimization problem reduces to solving the associated Fokker-Planck and its adjoint equations, the latter being non-trivial. We present a derivation procedure of the adjoint equation and its internal boundary conditions in time to efficiently derive the optimal transporting strategy. We finally provide a demonstrative computational example of a transporting problem of Ayu sweetfish (*Plecoglossus altivelis altivelis*) based on the latest real data set.

## 1. Introduction

### 1.1 Problem background

Inland fishery resources are important protein sources of human lives and their socioeconomic importance has been rapidly increasing [1-2]. Life histories of inland fishery resources, especially fishes migrating along rivers, are critically affected by transverse hydraulic structures like dams and weirs because they physically prevent the fishes from ascending/descending [3-4]. These physical barriers not only serve as obstructions of the migration processes, but also disturb fish assemblages in both their upstream and downstream reaches of rivers [5], resulting in fragmentation of fish habitats [6].

Environmental improvement schemes, such as barrier removal and flow adaptation, have been considered previously to mitigate the human-induced impacts [7-9]. Artificially releasing fishes from a habitat or an aquaculture farm to spatially apart target habitat(s) was considered as an alternative countermeasure, and has been applied to Chinook Salmon *Oncorhynchus tshawytscha* [10-11] and Ayu sweetfish *Plecoglossus altivelis altivelis* [12-13] so far. An advantage of this countermeasure is that it does not associate with water environmental regulations that are sometimes costly or even impossible. However, when and how much of the fishes should be transported is still an unresolved issue. To the best of the authors' knowledge, only a few reports for mountainous stream fishes exist [14]. This is the motivation of the current paper and we address the issue from a standpoint of optimization of partial differential equations (PDEs) considering the resource dynamics.

### 1.2 Mathematical background

Identifying growth dynamics of target fishery resources is of importance for optimizing their transporting strategies because individuals' size and/or population density influence their fisheries utilities [15-17]. Growth dynamics of individuals of fishery resources follow stochastic differential equations (SDEs) whose main governing variable is the size, such as body weight and body length. These SDEs are in general driven by continuous [18-19] and jump processes [20-21]. In any cases, the macroscopic appearance of individuals' body weights is considered as a (multiple of) a probability density function (PDF) governed by a Fokker-Planck equation [22] (FPE). FPEs are degenerate parabolic PDEs, implying that an optimization problem of transporting population governed by stochastic growth dynamics is formulated as a PDE-constrained optimization problem [23-24].

PDE-constrained optimization based on FPEs has been studied from the standpoint of maximum principle [25] where solving an optimization problem reduces to computing Fokker-Planck and adjoint equations concurrently. PDE-constrained optimization problems based on FPEs have various applications; they are including but are not limited to bilinear optimal control [26], model predictive control [27], cost-efficient switching [28], pedestrian motion control [29], equilibrium firm dynamics [30], and impulsive mean field games [31]. The aforementioned examples of the FPEs in PDE-constrained optimization have a theoretical idealization that the decision-maker, the controller of the target stochastic system, can intervene the system at any time; however, it is often difficult to continuously intervene the target dynamics especially when the dynamics are biological ones in natural environment [32-33]. Problems in inland fisheries are no exception. Recently, Yoshioka et al. [34] successfully overcome this difficulty by considered a cost-efficient transporting problem of fishery resource population based on the impulse (discontinuous) control where the decision-maker can control the target population dynamics impulsively without the need of using the conventional continuous control variables. This is an important engineering problem but has been neither mathematically nor numerically studied except in the above-mentioned recent research [34].

Another issue to be considered in applications is that modern growth models of individuals often have both observable and unobservable variables (variables that are difficult to directly observe). The decision-maker of the target system should make decisions based on observable variables. The unobservable variables are called latent variables in this paper. Examples of the observable variables are the body weight and body size, while the latent variables include the history-dependent maximum body weight [19, 35] and dynamically changing energetic variables [36-37]. Temporally inhomogeneous models can also be seen as open-ended growth models where time-dependent parameters add a flexibility in modeling the growth curves [38-41]. Hence, we need to find a mathematical framework that depends only on the observable variables. This is a new mathematical problem not approached so far.

### 1.3 Objective and contribution

Motivated by the research backgrounds described above, the objectives of this paper are to formulate and analyze a cost-efficient transporting problem of fishery resource populations based on impulsive control of FPEs. The stochastic growth model as a foundation of the optimization problem is inspired by an open-ended logistic model [19], where an individuals' growth curve is not based on a single SDE [18] but on a system of SDEs with which a more realistic probabilistic description of the stochastic body growth of fishes becomes possible. We present a simpler as well as more accurate model than the earlier model [19] containing simpler SDEs performing comparably well. A key ingredient of the new model is the novel use of the Wright-Fisher model famous in population genetics [42-43] as a latent driver of the growth dynamics, which turns out to generate statistically accurate stochastic growth curves fitting real data. This serves as a new application of the Wright-Fisher model to a completely different problem.

The corresponding FPE of the new model has smooth but degenerate coefficients. The control variables to stand for transporting strategies are set to be impulsive such that certain amount of the population in a habitat is transported to other habitat(s) at certain times. The decision-making is assumed to be based on either under a full-information case or a partial-information case; the latter is more realistic as it assumes decision-makings based on observable variables.

As in the existing framework of controlling FPEs [25], we solve the optimization problem using a maximum principle combined with an adjoint method. The impulse control assumption induces temporal interface conditions at prescribed discrete times [44-45], but the adjoint of the interface conditions is non-trivial even in the simplified case [34]. We therefore present a derivation procedure of the adjoint interface conditions and show that they are different between the full-information and partial-information cases. We show that solving the problem in the full-information case only requires solving the adjoint equation, while the partial-information case requires iteratively solving both Fokker-Planck and adjoint equations.

Impulsive controls are by virtue sparse in time and do not require continuously observing/intervening target dynamics. Such problems are not only interesting from a mathematical viewpoint, but also important in applications because they do not require continuous management efforts and are suited to the situation where the chances of interventions are severely limited. For example, occurrences of pandemics like the COVID-19 infections [46-47] critically limit the available labor forces for fisheries management both in quality and quantity.

We finally present an application example of the proposed model to *P. altivelis* as a major inland fishery resources in Japan. The optimization problems in full-information and partial-information cases are handled numerically using modern high-resolution solvers [48-49], where a Picard iteration is additionally utilized for handling the partial-information case to solve the coupled Fokker-Planck and adjoint equations. Consequently, we contribute to formulation, analysis, and application of a new PDE-based optimization problem under an assumption that not all the variables are observable. Our paper bridges scientific and engineering sides of the impulse control as a simple and efficient discontinuous control.

## 2. Mathematical model
### 2.1 Stochastic growth dynamics

We firstly present an existing model and its limitation. After that, we propose a new model that can overcome the limitation. The time is denoted as $t > 0$. The body weight of an individual fish at time $t$ is denoted as $X = (X_t)_{t \geq 0}$. This is a continuous-time and non-negative variable as biologically required. The model has another stochastic variable, which is the history-dependent maximum body weight $K = (K_t)_{t \geq 0}$ [19]. This is also a continuous-time and non-negative variable satisfying the almost sure (a.s.) relationships $X_t \leq K_t$ ($t \geq 0$) and $\lim_{t \to +\infty} X_t = \lim_{t \to +\infty} K_t$. Being different from the conventional logistic models where the maximum body weight is a fixed constant shared by the individuals [19], this model has an advantage that the different individuals possibly have different maximum body weights according to their life histories. In addition, the process $X$ should be increasing in time unless there are external environmental catastrophes that are not considered in this paper.

#### 2.1.1 Previous model

An open-ended system of Itô's SDEs governing the couple $(X, K)$ was proposed in Yoshioka et al. [19]:

$$\mathrm{d}X_t = r\left(1 - \frac{X_t}{K_t}\right)X_t \mathrm{d}t, \quad \mathrm{d}K_t = (X_t - K_t)(D\mathrm{d}t + \sigma \mathrm{d}B_t), \quad t > 0 \tag{1}$$

subject to an initial condition $(X_0, K_0)$ with $0 < X_0 < K_0$. Here, $r > 0$ is the specific growth rate, $\sigma > 0$ is the environmental noise intensity, $D > 0$ is a relaxation parameter, and $B = (B_t)_{t \geq 0}$ is the 1-D standard Brownian motion [50]. A deterministic prototype ($\sigma = 0$) is found in Thornley and France [35].

The term "open-ended" comes from the characteristic that sample paths of $(X_t)_{t \geq 0}$, even with the same initial condition, have different $\lim_{t\to+\infty} X_t$ in general. Interestingly, this natural property is not equipped with the conventional stochastic logistic model [18] having a common deterministic limit of $\lim_{t\to+\infty} X_t$. Empirically, the relaxation parameter $D$ modulates the speed of maturity: larger $D$ results in smaller limit $\lim_{t\to+\infty} X_t$. The noise intensity $\sigma$ represents stochastic fluctuation involved in the growth dynamics.

The system (1) admits a unique strong solution (path-wise continuous solution) (Theorem 2.1 of Yoshioka et al. [19]) and satisfies $\lim_{t\to+\infty} X_t K_t^{-1} = 1$. It reasonably captures the probability density of the body weight of *P. altivelis* collected in a river system during summer (growing season of the fish), while underestimates the mean growth in the coming autumn during which the fish spawns and dies out. See, also an example presented in **Section 4**.

We apply a transformation of variables so that the growth dynamics are described by the new couple $(W, Z)$, where $W = (W_t)_{t\geq 0}$ with $W_t = \ln X_t$ and $Z = (Z_t)_{t\geq 0}$ with $Z_t = X_t K_t^{-1}$. An advantage of using the transformed variables is that $Z_t$ is valued in the unit interval (compact set) $[0,1]$ being different from $K_t$ that is unbounded. Compact dynamics are easier to computationally implement. Another advantage is that the state-dependence of the drift coefficient of the SDE of the body weight becomes weaker than (1). By applying the classical Itô's formula [50] to (1), we get the equivalent system

$$dW_t = r(1-Z_t)dt \quad \text{and} \quad dZ_t = A(Z_t)dt + C(Z_t)dB_t, \quad t > 0 \qquad (2)$$

with the coefficients

$$A(Z_t) = Z_t(1-Z_t)(r + D + \sigma^2(1-Z_t)) \quad \text{and} \quad C(Z_t) = \sigma Z_t(1-Z_t), \qquad (3)$$

where the derivation procedure is presented in Section 2.2 of Yoshioka et al. [19]. The corresponding initial conditions are $W_0 = \ln X_0$ and $Z_0 = X_0 K_0^{-1}$.

### 2.1.2 Proposed model

A limitation of the existing model is a severe underestimation of the late growth of *P. altivelis*, our target species, as demonstrated in Yoshioka et al. [19]. The proposed model is simpler than the existing one, and further it does not face with the issue of underestimation (**Section 3.2**).

The proposed model directly solves SDEs of $(W, Z)$. As in the previous model, we formulate an SDE of $Z$ such that $Z_t$ ($t \geq 0$) is valued in the compact interval $[0,1]$. From an application viewpoint, a plausible model should be simulated without technical difficulties, especially computed sample paths of $Z$ should be in $[0,1]$ as well. It should also be able to capture growth dynamics of real data. We propose the following Itô's SDE model as a candidate complying with the requirements above:

$$dW_t = r(1-Z_t)dt \quad \text{and} \quad dZ_t = D(1-Z_t)dt + \sigma\sqrt{Z_t(1-Z_t)}dB_t, \quad t > 0 \qquad (4)$$

under the assumption $2D \geq \sigma^2$ and again $W_t = \ln X_t$. This system is equipped with an initial condition

$W_0 = \ln X_0$ such that $X_0 > 0$ and $Z_0 \in (0,1)$. The system (4) can be expressed in the form (2) with

$$A(Z_t) = D(1-Z_t) \text{ and } C(Z_t) = \sigma\sqrt{Z_t(1-Z_t)}. \tag{5}$$

The second equation of (4) is a Wright-Fisher model having a linear drift (Tran et al. [42], Chapter 6 of Alfonsi [43]) or more commonly the simplest case of polynomial jump-diffusions on the unit simplex [51].

Owing to the polynomial diffusion property and the fact that the dynamics of $Z$ is not affected by the dynamics of $X$, the system admits a unique strong solution $(Z_t)_{t \geq 0}$ bounded in $[0,1]$ globally in time (Theorem 6.1.1 of Alfonsi [43]). The first equation of (4) is a linear SDE having a uniformly bounded coefficient, which is trivially square integrable, admits a unique strong solution. The uniqueness of $(W_t)_{t \geq 0}$ directly means the uniqueness of $(X_t)_{t \geq 0}$ because of $X_t = \exp(W_t)$. The process $(W_t)_{t \geq 0}$ and thus $(X_t)_{t \geq 0}$ are increasing in time because of the boundedness property $Z_t \in [0,1]$ ($t \geq 0$). We can use numerical schemes to simulate bounded numerical solutions inheriting this property [52], further supporting usefulness of the proposed model in applications.

Due to the specific form of the drift coefficient $A$ in (5), applying a Feller's test to the second equation of (4) yields $\tau < +\infty$ a.s. (Proposition 6.1.2 of Alfonsi [43]), where $\tau$ is the first hitting time

$$\tau = \inf\{t \geq 0, Z_t = 1\}. \tag{6}$$

We have $\inf\{t \geq 0, Z_t = 0\} = +\infty$ owing to a McKean argument with the assumption $2D \geq \sigma^2$ (Exercise 6.1.3 of Alfonsi [43]), which is satisfied in a real case as demonstrated later.

Integrating the first equation of (4) yields the upper-bound satisfied at each $t \geq 0$:

$$W_0 \leq W_t = W_0 + \int_0^t r(1-Z_s)\mathrm{d}s \leq W_0 + rt. \tag{7}$$

Furthermore, the boundedness of $\tau$ implies $W_0 \leq W_\tau \leq W_0 + r\tau < +\infty$, meaning that growth of the fish in the proposed model terminates at a finite time with probability one.

In summary, we obtain the following well-posedness and boundedness results.

**Proposition 1** *The system (4) with an initial condition $(W_0, Z_0)$ such that $W_0 = \ln X_0$ with $X_0 > 0$ and $Z_0 \in (0,1)$ admits a unique strong solution $(W_t, Z_t)_{t \geq 0}$. This solution satisfies the following properties; $(W_t)_{t \geq 0}$ is increasing and $(Z_t)_{t \geq 0}$ is bounded in $[0,1]$ globally in time. Furthermore, this solution satisfies $\tau < +\infty$ with $\tau$ given in (6) as well as $\inf\{t \geq 0, Z_t = 0\} = +\infty$.*

**Remark 1** The assumption $2D \geq \sigma^2$, namely the assumption that the environmental noise intensity is sufficiently small, is necessary to have $\inf\{t \geq 0, Z_t = 0\} = +\infty$. In fact, $Z_s = 0$ at a time $s$ means $Z_t = 0$ ($t \geq s$), which leads to the convex and unbounded exponential growth $W_t = W_s + \int_s^t r\mathrm{d}u$ ($t \geq s$) that is critically different from biologically plausible sigmoid-like growth.

## 2.2 Fokker-Planck equation

The FPE governing the PDF of the system (2) is presented in this sub-section. The PDF of the variables $(t, W_t, Z_t) = (t, w, z)$ is expressed as $p = p(t, w, z)$. The computational domain of the FPE is truncated as $\Omega = (0, \bar{W}) \times (0, 1)$ where $\bar{W}$ is a sufficiently large value chosen so that the impact of truncating the domain is localized near the boundary $w = \bar{W}$. This localization technique is common in numerical computation of PDEs in unbounded domains (e.g., Jasso-Fuentes et al. [53]).

The governing FPE of $p$ is given by [22]

$$\frac{\partial p}{\partial t} + \frac{\partial}{\partial w} F_w(p) + \frac{\partial}{\partial z} F_z(p) = 0, \quad t > 0, \quad (w, z) \in \Omega \tag{8}$$

with the fluxes

$$F_w(p) = r(1-z)p \quad \text{and} \quad F_z(p) = A(z)p - \frac{\partial}{\partial z}\left(\frac{1}{2}C^2(z)p\right), \tag{9}$$

subject to the initial condition

$$p(0, w, z) = \delta(w - W_0, z - Z_0), \quad (w, z) \in \Omega, \tag{10}$$

where $\delta(w - W_0, z - Z_0)$ is the Dirac's Delta concentrated at $(W_0, Z_0) \in \Omega$. The couple $(A, C)$ of the coefficients can be given by (3) (the earlier model) or (5) (the proposed model), but we focus on the latter because the FPE of the former is not further discussed in this paper.

Boundary conditions must be specified along $w = 0, \bar{W}$ and $z = 0, 1$ when analyzing degenerate parabolic PDEs [54-55]. Because the domain is a rectangle, we can apply the Feller's condition [56] to each boundary. The zero-flux boundary condition $r(1-z)p = 0$ is specified along $w = 0$ to prevent entering the mass from outside the domain. The same boundary condition is prescribed along the other side $w = \bar{W}$ so that no mass goes out from this boundary. The zero-flux condition $A(z)p - \frac{\partial}{\partial z}\left(\frac{1}{2}C^2(z)p\right) = 0$ is prescribed along $z = 0, 1$ so that no probability mass goes out from this boundary. Actually, no boundary condition is needed along $z = 0$ where no sample paths reach by virtue of the SDEs. However, this boundary condition along $z = 0$ has been found to be innocuous if the mass is initially 0 along this boundary. The mass conservation $\int_\Omega p \, \mathrm{d}w \mathrm{d}z = 1$ ($t \geq 0$) is satisfied as well.

## 2.3 Population dynamics in each habitat

We continue to the population transportation problem based on FPEs. We do not directly use (8) but its variant considering mortality of the individuals as explained below. As in Yoshioka et al. [34], assume that there exist two habitats called $H_1$ and $H_2$, and that there is some population of a target fish species in each habitat. Assume further that these habitats are physically disconnected and that the decision-maker, a manager of the fishery resource, wants to transport the population from $H_1$ to $H_2$ due to some management reason. For example, this is a problem of releasing cultured fish to a natural river. Similar problems arise when considering transporting the population from a habitat with an inferior quality to that with a better quality. In what follows, the quantities related to the habitat $H_i$ ($i=1,2$) are denoted by the subscript $i$ like $r_i$. See also the conceptual diagram of the population dynamics (**Figure 1**).

By assuming population dynamics without density effects that can be justified under weak competition among individuals, we scale the total population $y_i$ (under no population transport and no reproduction during) by the corresponding PDF $p_i$ as $y_i = N_i e^{-R_i t} p_i$. Here, $R_i > 0$ is the mortality rate and $N_i > 0$ is the initial population at $t=0$. Theoretically, the population should be an integer variable but here we approximate it as a continuous variable because directly managing the discrete nature of the population is not practical from both theoretical and computational standpoints.

We consider transporting the population from $H_1$ to $H_2$. Assume that transportations are allowed only at prescribed times $\tau = \{\tau_j\}_{1 \le j \le J}$, $J \in \mathbb{N}$, $0 < \tau_1 < \tau_2 <,...,< \tau_J < T$, where $T > 0$ is a terminal time. At each $\tau_j$, the transporting strategy has the form [34]

$$\begin{pmatrix} y_1(\tau_j +) \\ y_2(\tau_j +) \end{pmatrix} = \begin{pmatrix} 1-u_j & 0 \\ u_j & 1 \end{pmatrix} \begin{pmatrix} y_1(\tau_j) \\ y_2(\tau_j) \end{pmatrix}, \quad 1 \le j \le J \tag{11}$$

with some control variable $u_j \in [0,U]$ with $0 < U < 1$. We used the abbreviations $y_i(\tau_j) = y_i(\tau_j, w, z)$ and $y_i(\tau_j +) = \lim_{t \to +0} y_i(\tau_j + t)$. The precise definition of $u = \{u_j\}_{1 \le j \le J}$ will be presented in the next sub-section. What is important here is that the transporting the population is formulated as a temporal interface condition (11).

Based on the discussion above, we obtain the population dynamics in each $H_i$ ($i=1,2$):

$$\frac{\partial y_i}{\partial t} + \frac{\partial}{\partial w} F_{w,i}(y_i) + \frac{\partial}{\partial z} F_{z,i}(y_i) = -R_i y_i, \quad 0 < t \notin \tau, \quad (w,z) \in \Omega \tag{12}$$

subject to the non-negative initial condition

$$y_i(0, w, z) = y_{i,0}(w, z), \quad (w,z) \in \Omega \tag{13}$$

with some integrable $y_{i,0}$ such that $\int_\Omega y_{i,0} \mathrm{d}w \mathrm{d}z = N_i \ge 0$ ($N_1 N_2 > 0$) and $y_{i,0} = 0$ on $w = 0$, coupled with the temporal interface condition (11). The equation (12) can be represented in an abstract form

$$\frac{\partial y_i}{\partial t} + S_i y_i = -R_i y_i \tag{14}$$

with the linear operator $S_i$ given by

$$S_i = \frac{\partial}{\partial w} F_{w,i}(\cdot) + \frac{\partial}{\partial z} F_{z,i}(\cdot). \tag{15}$$

Boundary conditions same with that explained in **Section 2.3** are specified as well. The PDE (12) is also called an FPE in what follows when there is no confusion. In fact, the equation (8) is not focused on below.

Recall that the variable $Z$ is latent. Therefore, the decision-maker can only observe $X$ in realistic cases. As the variable $y_i$ is the population density such that its integration over $\Omega$ gives the total population in $H_i$. For later use, set the conditional population density $\bar{y}_i = \bar{y}_i(t,w)$:

$$\bar{y}_i(t,w) = \int_0^1 y_i(t,w,z)\mathrm{d}z, \quad t > 0, \quad w \in \mathbb{R}. \tag{16}$$

This is the integrated population density with respect to the latent variable $z$, which turns out to be observable because the (logarithm of) body weight is assumed to be observable. The conditional population serves as a key quantity in the optimization problem under a partial-information assumption.

## 2.4 Objective function

The objective functions and the admissible sets of controls are presented for both full-information and partial-information cases.

### 2.4.1 Full-information case

The objective function contains the total cost of the transportation and the terminal reward of realizing the abundant population in the habitat $H_2$ in a preferred range $\omega \subset (0,+\infty)$ of the (logarithm of) body weight. Assume that the cost of transporting the population of the amount $u_j y_1$ is $cu_j y_1$ with a constant $c > 0$ [34]. The admissible set of controls is denoted as $\mathfrak{U}$ and is defined as

$$\mathfrak{U} = \left\{ u = \{u_j\}_{1 \le j \le J} \,\Big|\, u_j \in L^\infty(\Omega), 0 \le u_j \le U, 1 \le j \le J \right\}. \tag{17}$$

The objective function $\phi = \phi(y,u)$ is then set as

$$\phi(y,u) = \sum_{j=1}^{J} c \int_\Omega u_j(w,z) y_1(\tau_j) \mathrm{d}w\mathrm{d}z - \int_{\omega \times (0,1)} y_2(T) \mathrm{d}w\mathrm{d}z, \tag{18}$$

where the dependence of $y_i$ on the variables $(w,z)$ is omitted to simplify the description. Similar description will be used later if there will be no confusion. The first and second terms of $\phi$ represent the transportation cost and fishery utilities gained by achieving the preferred growth, respectively. The objective of the optimization problem is to choose an optimal control $u^* = \{u_j^*\}_{1 \le j \le J} \in \mathfrak{U}$ minimizing $\phi$.

***Remark 2*** The admissible set $\mathfrak{U}$ in (17) is bounded, closed, and convex. In addition, the objective function $\phi$ is bounded from below if all the integrals are bounded. Its minimizer, namely an optimal transporting strategy, exists under this assumption. This remark applies to the partial-information case below as well. See, also Theorem 4.1 of Mbogne and Thron [44] for a similar issue of a different model.

### 2.4.2 Partial-information case

As in the full-information case, the objective function is the sum of the total cost of the transportation and the terminal reward. The difference between the full-information and partial-information cases is that the admissible set of the controls is not (17) but

$$\mathfrak{U} = \left\{ u = \{u_j\}_{1 \leq j \leq J} \,\middle|\, u_j \in L^\infty(0, \overline{W}), 0 \leq u_j \leq U, 1 \leq j \leq J \right\}. \tag{19}$$

The corresponding objective function $\phi = \phi(y, u)$ is then set as

$$\begin{aligned}\phi(y,u) &= \sum_{j=1}^{J} c \int_{\Omega} u_j(w) y_1(\tau_j) \mathrm{d}w \mathrm{d}z - \int_{\omega \times (0,1)} y_2(T) \mathrm{d}w \mathrm{d}z \\ &= \sum_{j=1}^{J} c \int_0^{\overline{w}} u_j(w) \overline{y}_1(\tau_j) \mathrm{d}z - \int_\omega \overline{y}_2(T) \mathrm{d}z \end{aligned}. \tag{20}$$

The difference of the observability between the two cases appears in the objective functions as that of the partial-information case is described with the conditional population density $\overline{y}_i$, while that of the full-information case is not. This difference is inherited in the adjoint equations and the optimal controls.

***Remark 3*** One may want to use some governing PDE of the conditional population density instead of (12). This is because of a closure problem. In fact, we formally get the following equation directly from (12):

$$\frac{\partial \overline{y}_i}{\partial t} + R_i \overline{y}_i = -\frac{\partial}{\partial w} \int_0^1 F_{w,i}(y_i) \mathrm{d}z = -\frac{\partial}{\partial w}(r_i y_i) + \frac{\partial}{\partial w}\left(r_i \int_0^1 z y_i \mathrm{d}z\right). \tag{21}$$

The last term of (21) cannot be simplified further because of the need of the explicit information of $zy_i$ whose computation requires the information of $z^2 y_i$. We encounter an infinitely long cascade of higher-order moment problems that are not computable. One may overcome this issue using a data-driven closure approximation technique [57], but this is beyond the scope of this paper.

## 2.5 Adjoint equation and optimal control

The adjoint equation and the optimal control are presented for both full-information and partial-information cases. For the ease of explanation, the partial-information case is firstly explained.

### 2.5.1 Partial-information case

As a starting point to derive the adjoint equation and its temporal interface conditions, we calculate a directional derivative of the objective function $\phi$ with respect to the control variable using a sensitivity equation of the FPE. We show that the directional (Gâteaux) derivative of $\phi$ can be written compactly using the adjoint equation, and further that the optimal control is expressed using the adjoint variables. The directional differentiability of the population density $y_i$ is assumed below for simplicity of analysis. See, also **Section 2.6** on regularity of the FPE.

The sensitivity $s_i$ is defined as the directional derivative of $y_i$ at the control $u \in \mathfrak{U}$ in the direction of the small perturbation $v$ such that $v \in \mathfrak{U}$. The governing equation of $s_i$ is given by

$$\frac{\partial s_i}{\partial t} + S_i s_i = -R_i s_i, \quad 0 < t \notin \tau, \quad (w, z) \in \Omega \tag{22}$$

coupled with the interface condition

$$\begin{pmatrix} s_1(\tau_j +) \\ s_2(\tau_j +) \end{pmatrix} = \begin{pmatrix} 1-u_j & 0 \\ u_j & 1 \end{pmatrix} \begin{pmatrix} s_1(\tau_j) \\ s_2(\tau_j) \end{pmatrix} + \begin{pmatrix} -v_j & 0 \\ v_j & 0 \end{pmatrix} \begin{pmatrix} y_1(\tau_j) \\ y_2(\tau_j) \end{pmatrix}, \quad 1 \leq j \leq J \tag{23}$$

with the homogenous initial condition $s_i(0) \equiv 0$ and the same boundary condition with $y_i$ ($i = 1, 2$). Using $s_i$, the directional derivative $\delta\phi(y,u) \cdot v$ of $\phi$ is given by

$$\delta\phi(y,u) \cdot v = \sum_{j=1}^{J} c \left\{ \int_\Omega v_j(w) y_1(\tau_j) dwdz + \int_\Omega u_j(w) s_1(\tau_j) dwdz \right\} - \int_{\omega \times (0,1)} s_2(T) dwdz. \tag{24}$$

The sensitivity equation is equipped with the zero-flux condition analogous to that used in the FPE.

Now, we introduce an adjoint equation to rewrite the right-hand side of (24) so that the optimality condition is found efficiently. We also show that the temporal interface condition of the adjoint variables is related to the optimality condition. Set the formal adjoint $S_i^*$ of $S_i$ as

$$S_i^*(\cdot) = -r_i(1-z)\frac{\partial}{\partial w}(\cdot) - A_i(z)\frac{\partial}{\partial z}(\cdot) - \frac{1}{2}C_i^2(z)\frac{\partial^2}{\partial z}(\cdot). \tag{25}$$

The adjoint equation except at the intervention times is set as

$$-\frac{\partial q_i}{\partial t} + S_i^* q_i = -R_i q_i, \quad 0 < t \notin \tau, \quad (w, z) \in \Omega \tag{26}$$

with the terminal conditions $q_1(T) = 0$ and $q_2(T) = -\chi_{\omega \times (0,1)}$, where $\chi_{\omega \times (0,1)}$ is the indicator function of the set $\omega \times (0,1)$ such that it is equals to 1 if $(w,z) \in \omega \times (0,1)$ and to 0 otherwise.

The adjoint equation also requires proper boundary conditions. Again, following the Feller's condition [56], the homogenous Neuman condition $r_i(1-z)\frac{\partial q_i}{\partial w} = 0$ is applied to the boundaries $w = 0, \bar{W}$, while $\frac{1}{2}C_i^2(z)\frac{\partial q_i}{\partial z} = 0$ to the boundary $z = 1$. No boundary condition is necessary along $z = 0$. In the numerical computation and the analysis demonstrated later, we specify the homogenous boundary condition along all the boundaries because it has been found to be innocuous. The specific form of the adjoint equation turns out to be useful in efficiently solving the optimization problems as explained below. The following lemma on the adjoint relationship serves as a key identity of the rewriting procedure. In addition, it plays a fundamental role in determining interface conditions of the adjoint equation.

**Lemma 1** *We have the integral relationship*

$$\int_\Omega (S_i s_i) q_i dwdz = \int_\Omega s_i (S_i^* q_i) dwdz, \quad i = 1, 2, \quad 0 < t \notin \tau, \tag{27}$$

*provided that all the integrals are well-defined.*
**(Proof of Lemma 1)**

Owing to the prescribed boundary conditions, direct calculations show

$$\int_\Omega \frac{\partial}{\partial w} F_w(s_i) q_i \mathrm{d}w\mathrm{d}z = \int_\Omega \frac{\partial}{\partial w}\left(r_i(1-z)s_i\right) q_i \mathrm{d}w\mathrm{d}z = -\int_\Omega s_i r_i (1-z) \frac{\partial q_i}{\partial w}\mathrm{d}w\mathrm{d}z \tag{28}$$

and

$$\begin{aligned}
\int_\Omega \frac{\partial}{\partial z} F_z(s_i) q_i \mathrm{d}w\mathrm{d}z &= \int_\Omega \frac{\partial}{\partial z}\left(A_i(z)s_i - \frac{\partial}{\partial z}\left(\frac{1}{2}C_i^2(z)s_i\right)\right) q_i \mathrm{d}w\mathrm{d}z \\
&= -\int_\Omega s_i A_i(z) \frac{\partial q_i}{\partial z}\mathrm{d}w\mathrm{d}z + \int_\Omega \frac{\partial}{\partial z}\left(\frac{1}{2}C_i^2(z)s_i\right) \frac{\partial q_i}{\partial z}\mathrm{d}w\mathrm{d}z \\
&= -\int_\Omega s_i A_i(z) \frac{\partial q_i}{\partial z}\mathrm{d}w\mathrm{d}z - \int_\Omega s_i \frac{1}{2}C_i^2(z) \frac{\partial^2 q_i}{\partial z^2}\mathrm{d}w\mathrm{d}z
\end{aligned} \tag{29}$$

Combining (28) and (29) yields (27). The zero-flux condition for the FPE was used in (28), while both the zero-flux condition for the FPE and the Neumann condition for the adjoint equation were used in (29).

□

**Lemma 2** below gives a linkage between the sensitivity and adjoint variables at the intervention times.

*Lemma 2 We have the identity*

$$\sum_{j=1}^{J}\left(s_i(\tau_j+)q_i(\tau_j+) - s_i(\tau_j)q_i(\tau_j)\right) = s_i(T)q_i(T). \tag{30}$$

*(Proof of Lemma 2)*

For brevity, set $\tau_0 = 0$ and $\tau_{J+1} = T$. By an integration-by-parts in time, we obtain

$$\int_0^T\left(\int_\Omega s_i \frac{\partial q_i}{\partial t}\mathrm{d}w\mathrm{d}z\right)\mathrm{d}t + \int_0^T\left(\int_\Omega q_i \frac{\partial s_i}{\partial t}\mathrm{d}w\mathrm{d}z\right)\mathrm{d}t = \sum_{j=1}^{J+1}\int_{\tau_{j-1}}^{\tau_j}\left(\int_\Omega \frac{\partial(s_i q_i)}{\partial t}\mathrm{d}w\mathrm{d}z\right)\mathrm{d}t = \sum_{j=1}^{J+1}\left[s_i q_i\right]_{\tau_{j-1}}^{\tau_j}. \tag{31}$$

By **Lemma 1** and the sensitivity and adjoint equations except at the intervention times, the left-most side of (31) is rewritten as

$$\begin{aligned}
&\int_0^T\left(\int_\Omega s_i \frac{\partial q_i}{\partial t}\mathrm{d}w\mathrm{d}z\right)\mathrm{d}t + \int_0^T\left(\int_\Omega q_i \frac{\partial s_i}{\partial t}\mathrm{d}w\mathrm{d}z\right)\mathrm{d}t \\
&= \int_0^T\left(\int_\Omega s_i\left(R_i q_i + S_i^* q_i\right)\mathrm{d}w\mathrm{d}z\right)\mathrm{d}t + \int_0^T\left(\int_\Omega\left((-R_i s_i - S_i s_i)\right)q_i \mathrm{d}w\mathrm{d}z\right)\mathrm{d}t \\
&= 0
\end{aligned} \tag{32}$$

By (31)-(32), the initial condition of $s_i$, and the terminal condition of $q_i$, we obtain

$$\begin{aligned}
0 &= \sum_{j=1}^{J+1}\left[s_i q_i\right]_{\tau_{j-1}}^{\tau_j} \\
&= \sum_{j=1}^{J+1}\left(s_i(\tau_j)q_i(\tau_j) - s_i(\tau_{j-1}+)q_i(\tau_{j-1}+)\right) \\
&= s_i(\tau_{J+1})q_i(\tau_{J+1}) - s_i(\tau_0)q_i(\tau_0) - \sum_{j=1}^{J}\left(s_i(\tau_j+)q_i(\tau_j+) - s_i(\tau_j)q_i(\tau_j)\right) \\
&= s_i(T)q_i(T) - \sum_{j=1}^{J}\left(s_i(\tau_j+)q_i(\tau_j+) - s_i(\tau_j)q_i(\tau_j)\right)
\end{aligned} \tag{33}$$

leading to the desired result (30).

□

By **Lemma 2** and the definition of the adjoint variables, we obtain

$$\sum_{j=1}^{J}\left(s_1(\tau_j+)q_1(\tau_j+)-s_1(\tau_j)q_1(\tau_j)\right)=0 \tag{34}$$

and

$$\sum_{j=1}^{J}\left(s_2(\tau_j+)q_2(\tau_j+)-s_2(\tau_j)q_2(\tau_j)\right)=-s_2(T)\chi_{\omega\times(0,1)}. \tag{35}$$

Substituting (34) and (35) into (24) yields

$$\begin{aligned}\delta\phi(y,u)\cdot v &= \sum_{j=1}^{J}c\left\{\int_\Omega v_j(w)y_1(\tau_j)\mathrm{d}w\mathrm{d}z+\int_\Omega u_j(w)s_1(\tau_j)\mathrm{d}w\mathrm{d}z\right\}\\ &+\int_\Omega\sum_{j=1}^{J}\left(s_1(\tau_j+)q_1(\tau_j+)-s_1(\tau_j)q_1(\tau_j)\right)\mathrm{d}w\mathrm{d}z\\ &+\int_\Omega\sum_{j=1}^{J}\left(s_2(\tau_j+)q_2(\tau_j+)-s_2(\tau_j)q_2(\tau_j)\right)\mathrm{d}w\mathrm{d}z\\ &\equiv\sum_{j=1}^{J}\int_\Omega I_j\mathrm{d}w\mathrm{d}z\end{aligned} \tag{36}$$

where

$$\begin{aligned}I_j &= cv_j(w)y_1(\tau_j)+cu_j(w)s_1(\tau_j)\\ &+s_1(\tau_j+)q_1(\tau_j+)-s_1(\tau_j)q_1(\tau_j)+s_2(\tau_j+)q_2(\tau_j+)-s_2(\tau_j)q_2(\tau_j)\end{aligned}. \tag{37}$$

Substituting (23) into (37) yields

$$\begin{aligned}I_j &= v_j y_1(\tau_j)\{c-q_1(\tau_j+)+q_2(\tau_j+)\}\\ &+s_1(\tau_j)\{cu_j+(1-u_j)q_1(\tau_j+)-q_1(\tau_j)+u_j q_2(\tau_j+)\}\\ &+s_2(\tau_j)\{q_2(\tau_j+)-q_2(\tau_j)\}\end{aligned}. \tag{38}$$

This equality suggests specifying the temporal interface conditions of the adjoint equation as follows:

$$q_1(\tau_j)=q_1(\tau_j+)+u_j\left(c-q_1(\tau_j+)+q_2(\tau_j+)\right) \text{ and } q_2(\tau_j)=q_2(\tau_j+). \tag{39}$$

Here, recall that the adjoint equation must be integrated backward in time because it is a degenerate parabolic terminal and boundary value problem. By (39), we arrive at

$$I_j=v_j y_1(\tau_j)\{c-q_1(\tau_j+)+q_2(\tau_j+)\} \tag{40}$$

and thus

$$\delta\phi(y,u)\cdot v=\sum_{j=1}^{J}\int_\Omega v_j(w)y_1(\tau_j)\{c-q_1(\tau_j+)+q_2(\tau_j+)\}\mathrm{d}w\mathrm{d}z. \tag{41}$$

Now, we get the necessary optimality condition from (41). As in Theorem 3.2 of De los Reyes [24], we have the necessary optimality condition

$$\delta\phi(y,u)\cdot(v-u) = \sum_{j=1}^{J} \int_{\Omega} \{v_j(w) - u_j(w)\} y_1(\tau_j) \{c - q_1(\tau_j+) + q_2(\tau_j+)\} dwdz \geq 0 \quad (42)$$

for all $v \in \mathfrak{U}$. Assume that $y_1 \geq 0$, which is not a restrictive biological assumption. We have

$$\begin{aligned} &\int_{\Omega} \{v_j(w) - u_j(w)\} y_1(\tau_j) \{c - q_1(\tau_j+) + q_2(\tau_j+)\} dwdz \\ &= \int_0^{\bar{W}} \{v_j(w) - u_j(w)\} \int_0^1 y_1(\tau_j) \{c - q_1(\tau_j+) + q_2(\tau_j+)\} dwdz, \\ &\equiv \int_0^{\bar{W}} \{v_j(w) - u_j(w)\} L_j(w) dw \end{aligned} \quad (43)$$

meaning that, for each $w$, if $L_j(w) = \int_0^1 y_1(\tau_j) \{c - q_1(\tau_j+) + q_2(\tau_j+)\} dz \geq 0$ then $u_j^* = 0$. On the one hand, we have $v_j - u_j^* = v_j \geq 0$ in this case by $v_j \in [0, \bar{u}]$ and (42) is satisfied. On the other hand, if $L_j(w) < 0$, then we get $u_j^* = \bar{u}$. In fact, $v_j - u_j^* = v_j - \bar{u} \leq 0$ in this case by the requirement $v_j \equiv 0$ due to (42). Consequently, we see that the optimal control has the following feasible form:

$$u_j^* = u_j^*(w) = \begin{cases} 0 & (L_j(w) \geq 0) \\ \bar{u} & (\text{Otherwise}) \end{cases}. \quad (44)$$

This is justified if the set where $L_j(w) = 0$ has the Lebesgue measure zero. Namely, (44) is valid if the set of points such that $L_j(w) = 0$ has the width 0. Owing to (44), we derive the complete temporal interface condition of the adjoint equation as follows:

$$q_1(\tau_j) = q_1(\tau_j+) + u_j^*(w)(c - q_1(\tau_j+) + q_2(\tau_j+))u_j^*, \quad q_2(\tau_j) = q_2(\tau_j+). \quad (45)$$

Consequently, with the help of **Remark 2**, we obtain the following proposition to design the optimal transporting strategy between the two habitats. Notice that optimal transporting strategies may not be unique. In **Section 3.4**, we numerically demonstrate that the necessary optimality condition (42) leads to a reasonable transporting strategy.

***Proposition 2*** *An optimal transporting strategy minimizing the objective function in (18) is (44), which are obtained by solving the FPE (14) and adjoint equation (26).*

The FPE (14) and adjoint equation (26) are coupled through the interface condition (45), meaning that finding the optimal control $u^* = \{u_j^*\}_{1 \leq j \leq J}$ requires simultaneously solving the two equations. To tackle this issue, we use a naïve Picard iteration algorithm (**Algorithm 1**) that turns out to work effectively in numerical computation. Here, $K = 0, 1, 2, ...$ stands for the iteration number.

***Algorithm 1***

1. Set an initial guess $u_i^{(0)}$.
2. Set $K \to K + 1$.

3. Compute $y_i^{(K)}$.

4. Compute $q_i^{(K)}$ and $u_i^{(K)}$.

5. Compute the difference $\Delta^{(K)}$ between $u_i^{(K)}$ and $u_i^{(K-1)}$. If $\Delta^{(K)}$ is sufficiently small, then output $\left(y_i^{(K)}, q_i^{(K)}, u^{(K)}\right)$. Otherwise, go to the step 2.

The difference $\Delta^{(K)}$ in the step 5 is measured in some appropriate norm such as $l^1$ or $l^\infty$ norm. We use the $l^\infty$-norm in our numerical computation of the partial-information case presented in the next section. It completely vanishes when a candidate of a numerical solution converges because of the bang-bang nature of the optimal control in our problem. Convergence of the algorithm is achieved at a finite number of iterations, which are at most 5 in the computation below.

***Remark 4*** The sensitivity equation was used in the derivation procedure of the optimal control, but not necessary in the computation of the optimal control in **Algorithm 1**.

### 2.5.2 Full-information case

The sensitivity and adjoint equations in the full-information case is similar to that in the partial-information case except for the temporal interface condition. In the full-information case we arrive at

$$\delta\phi(y,u)\cdot v = \sum_{j=1}^{J}\int_{\Omega} v_j(w,z) y_1(\tau_j)\{c-q_1(\tau_j+)+q_2(\tau_j+)\}dwdz, \quad (46)$$

where each $v_j$ depends on both $w$ and $z$. This is qualitatively different from the partial-information case where each $v_j$ depends solely on $w$. The necessary optimality condition is then derived as

$$\delta\phi(y,u)\cdot(v-u) = \sum_{j=1}^{J}\int_{\Omega}\{v_j(w,z)-u_j(w,z)\} y_1(\tau_j)\{c-q_1(\tau_j+)+q_2(\tau_j+)\}dwdz \geq 0, \quad (47)$$

suggesting the optimal control that is feasible:

$$u_j^* = u_j^*(w,z) = \begin{cases} 0 & (c-q_1(\tau_j+)+q_2(\tau_j+)\geq 0) \\ \overline{u} & (\text{Otherwise}) \end{cases} \quad (48)$$

and the complete temporal interface condition of the adjoint equation

$$q_1(\tau_j) = q_1(\tau_j+)+u_j^*(w,z)\left(c-q_1(\tau_j+)+q_2(\tau_j+)\right)u_j^* \text{ and } q_2(\tau_j) = q_2(\tau_j+). \quad (49)$$

**Algorithm 1** is not necessary in the full-information case because the optimal control can be computed solely by the adjoint variables, as shown in (48).

## 2.6 Remarks on the optimization problems

Here, remarks on mathematical aspects of the optimization problems that we have not covered are presented. Our optimization problem consists of the interface conditions at the transporting times and the degenerate

PDEs between each consecutive transporting times. The interface conditions are well-defined if each $y_i$ is bounded and integrable in $\Omega$ at each time $t \geq 0$. The PDE parts need to be analyzed considering the fact that both the Fokker-Planck and adjoint equations are degenerate parabolic. For the FPE, due to the fully degenerate nature of the diffusion coefficient in the $w$ direction, the classical regularity estimates do not apply (e.g., Chapter 6 of Bogachev et al. [58]). Existence of distributional weak solutions belonging to the space of integrable functions $L^1\left(\left(\tau_{j-1}, \tau_j\right); \Omega\right)$ for each $i$, provided that $y_i\left(\tau_{j-1}+\right) \in L^1(\Omega)$ for each $j$, can be discussed following Theorem 3.1 of Dumont et al. [59] of a degenerate parabolic FPE. Well-posedness of the adjoint equations seem to be more subtle according to Remark 3.6 of Dumont et al. [59] for a degenerate PDE of a non-divergent form. Our problem is more complicated because the Fokker-Planck and adjoint equations are interacting, raising an issue to be addressed in future.

Fortunately, the coefficients of the PDEs are at most polynomials, which are smooth and Lipschitz continuous over $\Omega$, and can be smoothly extended to outside $\Omega$. One may regularize these equations by adding a small Laplacian term to these equations as in Mertz and Pironneau [60] to get smooth (and thus classical pointwise) solutions to the regularized counterparts, with which numerical computation of the optimization problem is expected to become easier. Adding a regularization term is a convenient way for easing the analysis of the optimization problem, justifying the directional differentiability of the population density, but introduces an error factor other than the domain truncation. Nevertheless, the diffusion coefficient for regularization can be arbitrary small if it is positive. Indeed, we show in the next section that numerical computation of the optimization problem works well without introducing such a regularization term.

## 3. Application
### 3.1 Target problem

We apply the proposed optimization problem to a realistic case to demonstrate its usability. The target problem is a planning a fish transportation scheme between upstream and downstream reaches of an existing dam of Hii River flowing in Shimane Prefecture, Japan. The river has the length of 153 (km) and the catchment area of 2,540 (km$^2$). The river starts from Sentsu-zan Mountain and finally pours to the Sea of Japan. Hii River has a multipurpose dam called Obara Dam having the associated reservoir called Sakura-Orochi Reservoir. Local fisheries cooperatives called Hii River Fisheries Cooperatives (HRFC) is authorizing fisheries in the midstream to upstream reaches of this river.

The fish *P. altivelis* is a major fishery resource of inland fisheries in Japan; Hii River is not an exception. Fish catches of the fish has been steadily decreasing during 2000s due not only to aging and population decrease of members of inland fisheries cooperatives, but also to degradation of river environment serving as habitats, which was triggered by the flow regulations by operating huge hydraulic structures like dams [61-62]. Especially, Obara Dam is completely fragmenting habitats of aquatic species including *P. altivelis* because of its heigh 90 (m) as well as the lack of fish ladders.

Currently, designing effective release policy of juveniles from a training facility or from some river to a target river has been a central issue to sustain wild population and fishery catch of *P. altivelis* in Japan. Each release event can be impulsive because it takes few hours to one day in most cases. This applies to annual release projects conducted in Hii River and has also been adopted in the conceptual model for designing release policy of the fish [63]. In the framework of our model, the former can be set as $H_1$ while the latter as $H_2$.

The population of *P. altivelis* in Hii River has long been considered to be supported by the released populations in every year; however, the authors and HRFC found that the land-locked *P. altivelis* having a unique life-cycle exists in the upstream river of Sakura Orochi Reservoir. Detailed population dynamics of the land-locked *P. altivelis* is currently under investigation. At the same time, HRFC is planning to catch juveniles of the land-locked fish in the reservoir during spring and distribute them to other places of the same river such as river reaches downstream of the dam. There exists no report of density-dependent growth of *P. altivelis* in Hii River, motivating us to apply the proposed impulse control model for designing release policy of the fish in this river.

### 3.2 Parameter identification of the growth model

We use model parameter values based on an identified model of the fish growth at Hii River in 2019. We have the two data set. The first data set is a histogram (discrete PDF), which is called the histogram data in this paper, collected at August 4 in 2019 by union members of HRFC. The second data set is a time series data having daily-averaged body weights of the fish from July to November in 2019 collected by a member of HRFC, which is called the historical data in this paper. The sampling period of the historical data is irregular because he/she is a part-time fisher. In addition, the historical data has only non-chronological data [64] because each caught fish is not released back to the river. Therefore, it is not straightforward to identify the parameter values using the historical data. Consequently, we firstly identify the parameter values of the growth dynamics using the histogram data, and then verify the identified model against the historical data. In this way, we identify a reasonable model based both on the histogram and historical data.

As explained above, the model parameters of the system (4) is identified using the histogram data. In total 227 individuals were caught on August 4 in 2019. The maximum and minimum values of the body weight among the individuals are 20.0 (g) and 119.5 (g), respectively. In this paper, the model parameters have been identified so that the following performance measure $P$ based on the average (Ave), standard deviation (Std), and skewness (Skw) is minimized:

$$P = \frac{\left|\text{Ave}_{\text{Obs}} - \text{Ave}_{\text{Model}}\right|}{\text{Ave}_{\text{Obs}}} + \frac{\left|\text{Std}_{\text{Obs}} - \text{Std}_{\text{Model}}\right|}{\text{Std}_{\text{Obs}}} + \frac{\left|\text{Skw}_{\text{Obs}} - \text{Skw}_{\text{Model}}\right|}{\text{Skw}_{\text{Obs}}}. \tag{50}$$

Here, the sub-scripts Obs and Model represent the statistics computed using the observed data and that computed using the model. This performance measure is a sum of the relative errors of the first-order to the third-order statistics and a minimizing model should be able to capture the distributional shape of the observed histogram both qualitatively and quantitatively. The system (4) is numerically simulated using a Monte-Carlo method based on the bounded scheme [52] for the SDE of $Z$ and the classical Euler-Maruyama scheme for the SDE of $W$. Each sample path generated by the Monte-Carlo method then rigorously satisfies the bound $0 \leq W_t \leq 1$ without artificially truncating numerical solutions. The total number of sample paths is 1,000,000 and the time increment for discretization of the system is 0.004 (day). Increasing the computational resolution does not critically affect our identification result. The day of observation (August 4 in 2019) is set to be 90 (day) because of the assumption that the *P. altivelis* is released in the river around the beginning of May with the mean body weight of 6 (g). We thus set $X_0 = 6$ (g).

**Figure 2** presents the histogram of the observed data and that generated by an manually identified model whose parameter values, which are $r = 0.051$ (1/day), $D = 0.019$ (1/day), $\sigma = 0.051$ (1/day$^{1/2}$) $X_0 = 6$ (g), $W_0 = 0.02$ (-), suggesting that the proposed model can capture the probability law of the growth of *P. altivelis* in Hii River at the summer day in 2019. Especially, the unimodal and skewed nature of the observed histogram is correctly reproduced with the proposed model (**Table 1**). The identification results suggest that the model identified by optimizing the error index (50) gives a model that is consistent with the observation data of the growth of the fish *P. altivelis* in Hii River.

The identified model is then examined against the historical data as shown in **Figure 3**. The total number of the historical data in 2019 is 57. The comparison results show that the proposed model captures the mean behavior of the observed historical data. Recall that this model was identified solely with the histogram data, suggesting that the identified model is reasonable one for analyzing growth dynamics of *P. altivelis* in Hii River. A computational result of the earlier model with the optimized parameter values (Chapter 2 of Yoshioka et al. [65]) is also plotted in **Figure 3**, clearly showing that it underestimates the late growth period of the fish.

Because of considering a transportation problem of fish populations from one habitat to another, we also prepare another set of parameters to be used in our demonstrative numerical computation. Field surveys to collect the body weight data of *P. altivelis* were also planned in 2020, but the survey to collect histogram data during summer period was not cancelled due to the COVID-19 outbreak. Fortunately, a survey to collect historical data from July to November in 2020 was completed safely. We specify the growth model in 2020 using this historical data. We assume that the difference between the growth dynamics between 2019 and 2020 is generated by the intrinsic growth rate $r$, and explored the value of $r$ minimizing the least-squares error between the computed mean $\mathbb{E}[X_t]$ by the Monte-Carlo method and the observed body weights. We set the following minimization problem to identify $r$ in 2020:

$$r = \arg\min\left\{\frac{1}{N}\sum_{i=1}^{N}\left(\mathbb{E}\left[X_{t_i}\right] - \bar{X}_{t_i}\right)^2\right\} \equiv \arg\min\{\text{Err}\}, \qquad (51)$$

where $N = 73$ is the total number of the historical data collected in 2020.

We have identified the intrinsic growth rate in 2020 as $r = 0.048$ (1/day) as indicated in **Table 2**, showing that the growth rate is smaller in 2020 than 2019. **Figure 4** compares the statistics of identified models and the observation results, showing that the observation results are scattered around the mean curve of the identified model and almost falling between the curves of Average ± Standard deviation.

### 3.2.1 Other model parameters

The domain $\Omega$ is set as $(0, \overline{W}) \times (0,1)$ with $\overline{W} = 5.3$ corresponding to the maximum body weight 200 (g), which is a sufficiently large value compared with the observed body weights of *P. altivelis* in Hii River so far. Set the terminal time as $T = 70$ (day) and the intermediate times at which transporting the population is allowed as $\tau_j = 10j$ (day) ($j = 1, 2, 3, 4, 5, 6$). The target set $\omega$ is set as $(0.3\overline{W}, 0.7\overline{W})$. The initial condition of the FPE is set as follows:

$$y_1(0, w, z) = \frac{N}{\overline{W}} \exp\left[-a\left\{\left(\frac{w - \overline{w}}{\overline{W}}\right)^2 + (z - \overline{z})^2\right\}\right] \quad \text{and} \quad y_2(0, w, z) = 0, \quad (w, z) \in \Omega \tag{52}$$

with the parameter values $a = 500$, $\overline{w} = 0.1$, $\overline{z} = 0.1$, and $N$ is determined so that integrating $y_1(0, w, z)$ over the domain $\Omega$ becomes $6 \times 10^6$ (g). This assumes that, at the initial time $t = 0$, the individual's weights are $O(10^0)$ (g) and the total number of the individuals in the habitat $H_1$ is $O(10^6)$. However, because of the linearity of the objective function on the population, we get the same optimal control if we use other values of $N > 0$. Using the regularized initial condition (52) instead of a delta function is to avoid possible regularity issues of the FPE caused by its degenerate coefficients. For the transporting strategy, we set $U = 0.2$.

The parameter $c$ is firstly set as 0.2, but different values of $c$ are examined later. The parameters for the population dynamics are set as follows. We set $D_1 = D_2 = 0.019$ (1/day) and $\sigma_1 = \sigma_2 = 0.051$ (1/day$^{1/2}$) based on the identification results presented above. The mortalities are set as $R_1 = R_2 = 0.01$ (1/day). By contrast, we assume the different values of the intrinsic growth rates $r$ (1/day) between the habitats as follows: $(r_1, r_2) = (0.051, 0.048)$ or $(0.048, 0.051)$ (1/day).

### 3.3 Numerical method and computational conditions

We numerically discretize the Fokker-Planck and adjoint equations using the local Lax-Friedrichs schemes equipped with the conservative Weighted Essentially Non-Oscillatory (WENO) reconstruction [48] and the non-conservative WENO reconstruction [49], respectively. The temporal discretization is carried out using a classical Heun method (also called second-order Runge-Kutta method) and the diffusion terms are discretized using the classical second-order central difference method. Similar numerical methods have already been successfully applied to linear and nonlinear test cases extensively, and recently to a single-variable full-information case [34]. The WENO reconstructions accurately manage the first-order spatial partial differential terms that would be the largest source of errors in computing degenerate parabolic PDEs. Of course, one can use any other discretization methods similarly, such as the pseudo-spectral method [66] and the exponential difference method [67] if preferred.

WENO discretization can archive high-resolution but is non-monotone, meaning that they generate numerical solutions having spurious oscillations although they have been found to be small. These oscillations are harmful if numerical solutions must be rigorously bounded in a range or if one's interest is in computing solutions having monotone profiles [68-69]. In our case, the relationship (45) is based on the non-negativity of the population density that is not necessarily satisfied in numerical solutions computed using the conventional WENO discretization. We equip the WENO reconstruction for the FPE with a bound-preserving limiter [70], so that the computed population density provably stays non-negative. This limiting technique has been employed for high-resolution numerical computation of engineering problems based on PDEs [71-73] where physical constraints of the state variables such as positivity or non-negativity of the density must be numerically satisfied for stable computation.

We are therefore using an existing numerical scheme to the new optimization problem proposed in this paper. The time increment is set as $0.01$ (day) and the domain $\Omega$ is discretized by $201 \times 201$ uniform mesh that has been preliminary found to be sufficiently fine.

### 3.4 Computational results

We firstly compare the optimally transported populations between the full-information and partial-information cases; the former is more flexible but less realistic. Secondly, we compare the optimal transporting strategies of the partial-information model for different values of the transportation cost.

We firstly consider the case where the population is transported from a bad habitat $H_1$ to a good habitat $H_2$. Hence, we choose $(r_1, r_2) = (0.048, 0.051)$ (1/day). **Figures 5-6** show the conditional populations $\bar{y}_i = \bar{y}_i(t, w)$ $(i = 1, 2)$ for the full-information case. **Figure 7-8** show the adjoint variables $q_i = q_i(t, w, z)$ $(i = 1, 2)$ and the corresponding optimal transporting strategies at the time $t = \tau_3 = 30$ (day). Sharp transitions of the adjoint variables are captured without visible spurious oscillations. The computational results indicate that the population in the habitat $H_2$ is concentrated in the targeted range $\omega = (0.3\bar{W}, 0.7\bar{W})$ near the terminal time $T = 70$ (day) by the control depending both on $(W, Z)$. However, recall that $Z$ is latent and thus is difficult or even impossible to observe.

Similarly, **Figures 9-10** show the conditional populations $\bar{y}_i = \bar{y}_i(t, w)$ $(i = 1, 2)$ with the corresponding optimal transporting strategies for the partial-information case. The computational results show that the population in the habitat $H_2$ is concentrated in the targeted range $\omega$ near the terminal time $T$ in this case as well, although the control is now dependent only on the observable variable $W = \ln X$. The undulating profiles of the populations are due to the impulsive transportations of the population between the two habitats. In addition, the **Figures 9-10** imply that the optimal transporting strategy at each time $\tau_j$ $(j = 1, 2, 3, 4, 5, 6)$ is either an interval or null. **Figures 6** and **10** show that the transported populations have close profiles between the full- and partial- information cases. Hence, disadvantage of having the latent variable sems not to be significant in case of *P. altivelis* in this case.

We get qualitatively similar results for the parameterization representing $(r_1, r_2) = (0.051, 0.048)$ (1/day) where the population should be transported from a good habitat $H_1$ to a bad habitat $H_2$. Such a case would arise if the fishing is not allowed in $H_1$ but in $H_2$ according to fisheries regulations. As an example, of *P. altivelis* in Hii River, there are several areas closed to fishing that have dangerous terrain but high food availability for the fish.

We further analyze optimal transporting strategies in the partial-information case, especially its dependence on the transportation cost $c$. **Figures 11-12** show the optimal transporting strategies for different values of $c = 0.1k$ $(k = 1, 2, ..., 10)$ for the parameterizations $(r_1, r_2) = (0.048, 0.051)$ and $(r_1, r_2) = (0.051, 0.048)$, respectively. The computational results for $c = 0.9$ and $c = 1.0$ are not appearing in the figures because they are completely null in both parameterizations of $(r_1, r_2)$, suggesting that the fish population should not be transported if the transportation cost is huge. The interval where transporting the population is optimal shrinks as $c$ increases and eventually vanishes. This tendency is more significant for smaller $j$, implying that transporting the population should focus on the later chances closer to the terminal time $t = T = 70$ (day) as the cost increases. By comparing **Figures 11** and **12**, we find that smaller individuals should be transported with $(r_1, r_2) = (0.051, 0.048)$ than that with $(r_1, r_2) = (0.048, 0.051)$ because of the less possibility of larger growth of the transported individuals in the former case. This suggests a recommendation that efficiently transporting the fish population should consider the growth rate in both habitats.

As demonstrated in this section, depending on the estimated growth dynamics model of a fishery resource and an objective function considering the cost and benefit of transporting them, the optimal transporting strategy can be numerically computed by solving the Fokker-Planck and adjoint equations.

## 4. Conclusion

We formulated and analyzed a new PDE-constrained problem of impulsively transporting fishery resources based on the degenerate parabolic 2-D FPE and the associated adjoint equation. Its key elements were the new open-ended stochastic growth model having a latent driver and an efficient impulsive control formulation. The full-information and partial-information cases were discussed and computed focusing on the recent *P. altivelis* management case in Japan. The proposed growth model fitted the real data well, and the optimal transporting strategies were successfully obtained using the high-resolution schemes. The computational results suggested adaptive transporting strategies depending both on time and individuals' body weights. Mathematical modeling and computational analysis carried out in this paper will be able to provide useful information about cost-efficient transporting strategies of fishery resources.

We only analyzed a population transport problem between two habitats, but the proposed mathematical framework can be applied to more general problems transporting populations from habitat(s) to other habitat(s), such as an allocation problem of cultured fish to multiple points in a river system [13]. Such problems can be formulated based on the proposed framework but require higher computational costs and more intensive field surveys for data collection. The doubly stochastic case where the macroscopic population dynamics are perturbed by uncertain environmental disturbance can also be a critical issue arising in real problems. Considering a stochastic-PDE based optimization problems [74-75] will be a viable way to approach this issue. Regularity issues of the Fokker-Planck and adjoint equations that were not addressed in this paper are key future topics as well. Optimization models with adaptively choosing the opportunities to transport fishery resources can also be developed based on the proposed model. Identification of a more realistic initial condition of the population density and analysis of the degenerate parabolic and impulsively-controlled FPEs subject to it will also be key issues.

**Conflict of Interest Statement**
The authors do not have competing of interests to declare.

**Tables**

**Table 1.** Comparison of the observed and modelled statistics.

| Statistics | Observation | Model | Relative error |
|---|---|---|---|
| Average (g) | 56.4 | 56.4 | 0.000 |
| Standard deviation (g) | 18.2 | 18.3 | 0.005 |
| Skewness (-) | 0.95 | 0.94 | 0.011 |

**Table 2.** Comparison of Err to identify $r$ (1/day) in 2020. The underline corresponds to the identified model in 2020.

| $r$ | Err |
|---|---|
| 0.041 | 580.5 |
| 0.042 | 479.6 |
| 0.043 | 386.7 |
| 0.044 | 303.6 |
| 0.045 | 232.7 |
| 0.046 | 176.6 |
| 0.047 | 138.5 |
| <u>0.048</u> | <u>121.9</u> |
| 0.049 | 130.7 |
| 0.050 | 169.5 |
| 0.051 | 243.7 |

**Figures**

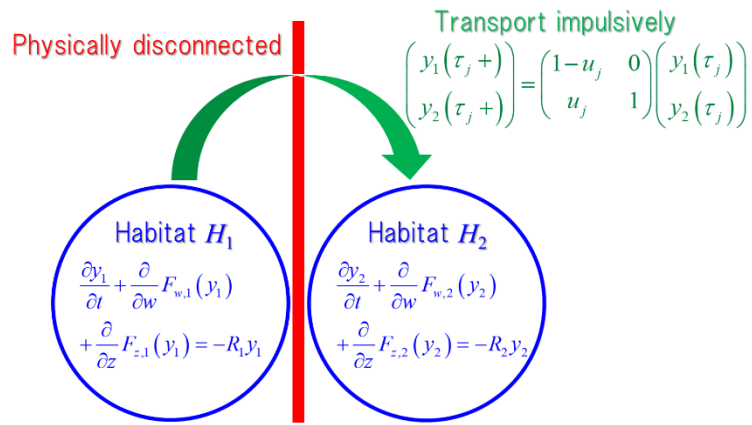

**Figure 1.** Conceptual diagram of the population dynamics.

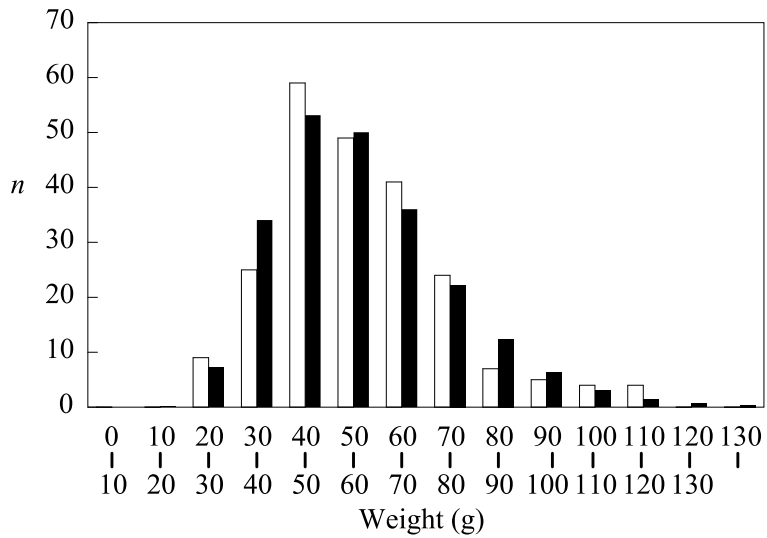

**Figure 2.** Comparison of the observed (White bars) and modelled histograms (Black bars). The label "0 - 10" represents the range $[0,10)$ (g). The same rule applies to the other labels except for the right-most one representing the unbounded range $[130,+\infty)$ (g).

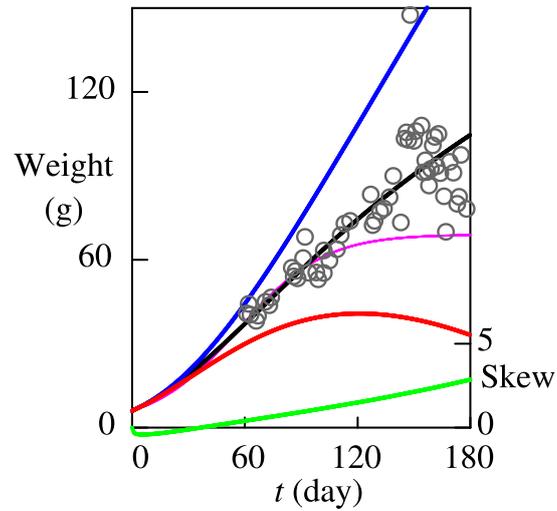

**Figure 3.** Comparison of the computed (Curves) and observed (Circles) historical data in 2019: Average (Black), Average + Standard deviation (Blue), Average − Standard deviation (Red), and Skewness (Green). The violet curve represents Average of the earlier model with the manually optimized parameter values ($r = 0.038$ (1/day), $D = 0.016$ (1/day), $\sigma = 0.099$ (1/day$^{1/2}$), $X_0 = 6$ (g), $W_0 = 0.03$ (Chapter 2 of Yoshioka et al. [66]), demonstrating that it underestimates the mean growth in this year.

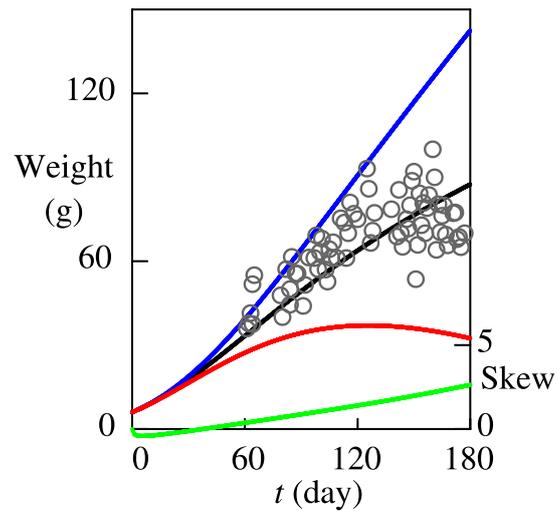

**Figure 4.** Comparison of the computed (Curves) and observed (Circles) historical data in 2020: Average (Black), Average + Standard deviation (Blue), Average − Standard deviation (Red), and Skewness (Green).

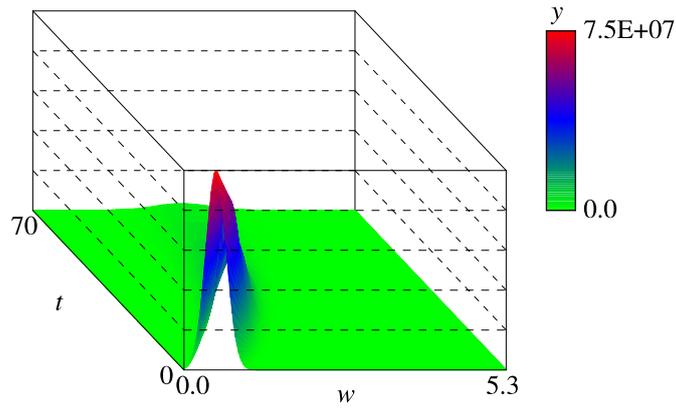

**Figure 5.** Computed conditional population $\bar{y}_1 = \bar{y}_1(t,w)$ for the full-information case.

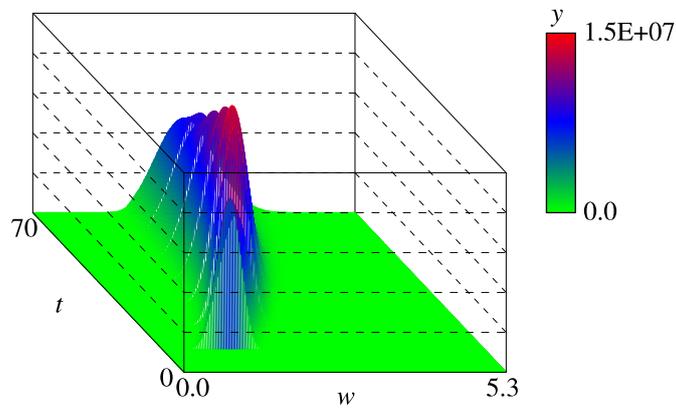

**Figure 6.** Computed conditional population $\bar{y}_2 = \bar{y}_2(t,w)$ for the full-information case.

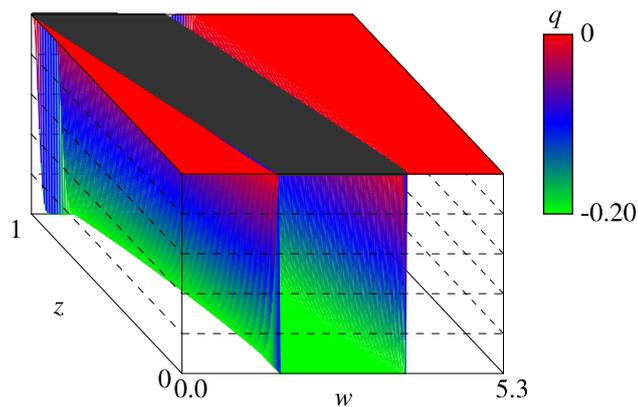

**Figure 7.** Computed adjoint variable $q_1 = q_1(t,w,z)$ (colored surface) and the corresponding optimal control (circle plots represent the area where the transportation is activated) for the full-information case.

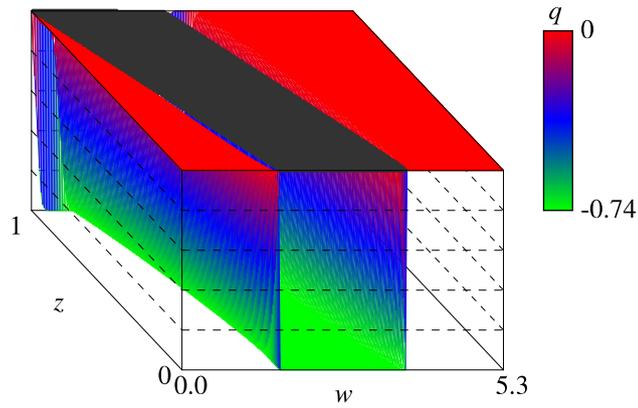

**Figure 8.** Computed adjoint variable $q_2 = q_2(t,w,z)$ and the corresponding optimal control for the full-information case. The same legends with **Figure 6**.

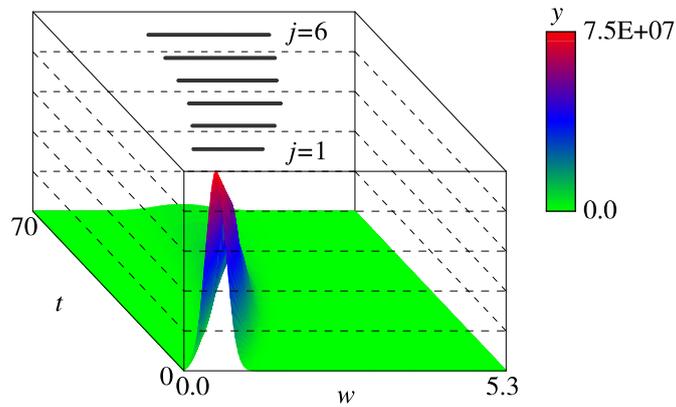

**Figure 9.** Computed conditional population $\bar{y}_1 = \bar{y}_1(t,w)$ (colored surface) with the corresponding optimal transporting strategy (grey plots) for the partial-information case.

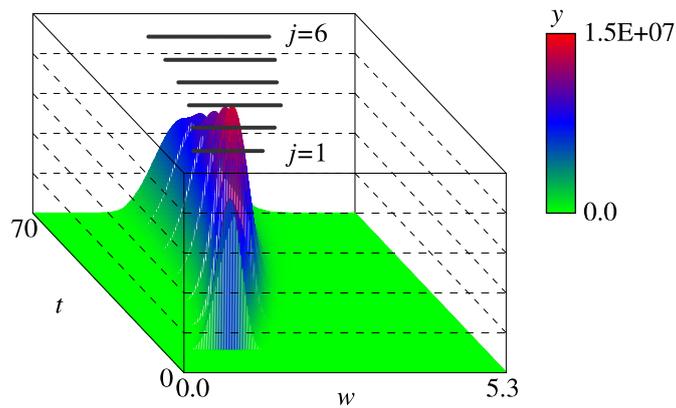

**Figure 10.** Computed conditional population $\bar{y}_2 = \bar{y}_2(t,w)$ with the corresponding optimal transporting strategy for the partial-information case. The same legends with **Figure 8**.

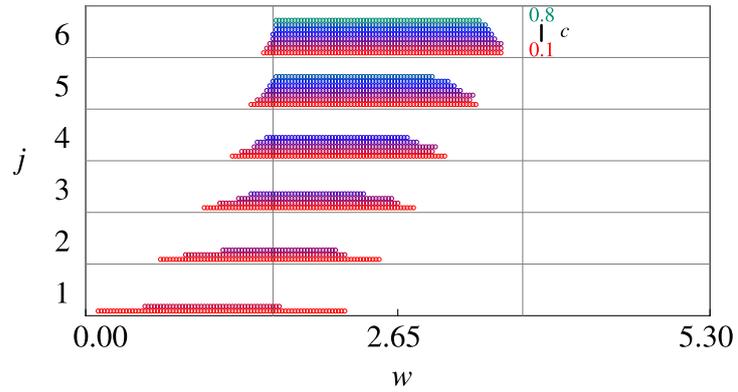

**Figure 11.** Computed optimal transporting strategies for different values of $c = 0.1k$ $(k = 1, 2, ..., 10)$ with $(r_1, r_2) = (0.048, 0.051)$ (1/day).

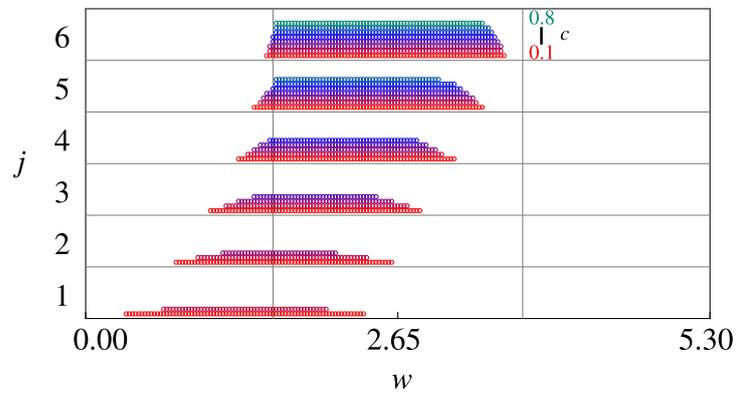

**Figure 12.** Computed optimal transporting strategies for different values of $c = 0.1k$ $(k = 1, 2, ..., 10)$ with $(r_1, r_2) = (0.051, 0.048)$ (1/day).